\documentclass[11pt,a4paper]{article}
\usepackage{amsmath, amsthm, amsfonts, amssymb, color}
\usepackage{mathrsfs}
\usepackage{color}
\usepackage{bbm}
\usepackage[notcite,notref]{}

\theoremstyle{definition}

\newcommand{\scr}[1]{\mathscr #1}
\definecolor{wco}{rgb}{0.5,0.2,0.3}

\numberwithin{equation}{section} \theoremstyle{remark}

\newcommand{\ua}{\uparrow}

\title{{\bf Bismut Formula and Gradient Estimates for Dirichlet Semigroups with  Application to Singular Killed DDSDEs }
}
\author{
	{\bf   Feng-Yu Wang, Xiao-Yu Zhao  }\\
	\footnotesize{Center for Applied Mathematics and KL-AAGDM, Tianjin University,   300072, China}\\
	\footnotesize{ wangfy@tju.edu.cn, zhxy\_0628@tju.edu.cn}}
\begin{document}
	\allowdisplaybreaks
	\def\R{\mathbb R}  \def\ff{\frac} \def\ss{\sqrt} \def\B{\mathbf
		B}\def\TO{\mathbb T}
	\def\I{\mathbb I_{\pp M}}\def\p<{\preceq}
	\def\N{\mathbb N} \def\kk{\kappa} \def\m{{\bf m}}
	\def\ee{\varepsilon}\def\ddd{D^*}
	\def\dd{\delta} \def\DD{\Delta} \def\vv{\varepsilon} \def\rr{\rho}
	\def\<{\langle} \def\>{\rangle} \def\GG{\Gamma} \def\gg{\gamma}
	\def\nn{\nabla} \def\pp{\partial} \def\E{\mathbb E}
	\def\d{\text{\rm{d}}} \def\bb{\beta} \def\aa{\alpha} \def\D{\scr D}
	\def\si{\sigma} \def\ess{\text{\rm{ess}}}
	\def\beg{\begin} \def\beq{\begin{equation}}  \def\eed{\end{equation}}\def\F{\scr F}
	\def\Ric{{\rm Ric}} \def\Hess{\text{\rm{Hess}}}
	\def\e{\text{\rm{e}}} \def\ua{\underline a} \def\OO{\Omega}  \def\oo{\omega}
	\def\tt{\tilde}
	\def\cut{\text{\rm{cut}}} \def\P{\mathbb P} \def\ifn{I_n(f^{\bigotimes n})}
	\def\C{\scr C}      \def\aaa{\mathbf{r}}     \def\r{r}
	\def\gap{\text{\rm{gap}}} \def\prr{\pi_{{\bf m},\varrho}}  \def\r{\mathbf r}
	\def\Z{\mathbb Z} \def\vrr{\varrho} \def\ll{\lambda}
	\def\L{\scr L}\def\Tt{\tt} \def\TT{\tt}\def\II{\mathbb I}
	\def\i{{\rm in}}\def\Sect{{\rm Sect}}  \def\H{\mathbb H}
	\def\M{\scr M}\def\Q{\mathbb Q} \def\texto{\text{o}} \def\LL{\Lambda}
	\def\Rank{{\rm Rank}} \def\B{\scr B} \def\i{{\rm i}} \def\HR{\hat{\R}^d}
	\def\to{\rightarrow}\def\l{k_1}\def\iint{\int}
	\def\EE{\scr E}\def\Cut{{\rm Cut}}\def\W{\mathbb W}
	\def\A{\scr A} \def\Lip{{\rm Lip}}\def\S{\mathbb S}
	\def\BB{\mathbb B}\def\Ent{{\rm Ent}} \def\i{{\rm i}}\def\itparallel{{\it\parallel}}
	\def\g{{\mathbf g}}\def\Sect{{\mathcal Sec}}\def\T{\mathcal T}\def\V{{\bf V}}
	\def\PP{{\bf P}}\def\HL{{\bf L}}\def\Id{{\rm Id}}\def\f{{\bf f}}\def\cut{{\rm cut}}
	\def\Ss{\mathbb S}
	\def\BL{\scr A}\def\Pp{\mathbb P}\def\Pp{\mathbb P} \def\Ee{\mathbb E} \def\8{\infty}\def\1{{\bf 1}}
	\maketitle
	
	\begin{abstract}
		By establishing a local version of  Bismut formula for Dirichlet semigroups on a regular domain, gradient estimates are derived for killed SDEs with  singular drifts. 	As   an  application, the total variation distance
		between two solutions of killed DDSDEs is bounded above  by  the  truncated $1$-Wasserstein distance  of initial distributions, in the regular and singular cases respectively.  		
	\end{abstract}
	
	\section{Introduction  }
	
	The Bismut formula, first established  in	Bismut \cite{Bismut} for the heat semigroup on Riemannian manifolds by using Malliavin calculus, has been intensively studied and applied for   SDEs and SPDEs.
	In order to derive gradient estimates on diffusion semigroups by using local geometry quantities, a local version of  Bismut formula is established in \cite{A} by refining the martingale argument developed in \cite{EL}, which is then  applied in \cite{AW} to derive explicit gradient estimates on Dirichlet semigroups.
	
	In recent years, the Bismut formula has been established in \cite{XXZZ} for singular SDEs where the drifts only satisfy a local integrability condition.
	When the drifts  also depend on the distribution of the solution, the SDEs are known as distribution dependent SDEs (DDSDEs for short). In this case,   the Bismut formula has been established in \cite{BRW20, HSW21, RW19, W23b} for the intrinsic/Lions derivative with respect to the initial distribution, as well as in \cite{R25} for the extrinsic derivative.
	
	To characterize nonlinear Dirichlet problem, the well-posedness and Lipschitz   continuity  in initial distribution have been derived in \cite{W23a}  for a class of singular DDSDEs   in a domain with killed boundary. As far as we know,  there is no any results on Bismut formula and applications for   killed distribution dependent SDEs. In this paper, we intend to establish the Bismut formula for Dirichlet semigroups for killed SDEs,
	and apply it to derive gradient estimates on  DDSDEs.
	
	Before moving on, let us recall the gradient estimate derived in \cite{AW} for the Dirichlet semigroup on manifold, and formulate the equivalent estimate on  distributions of the associated killed diffusions.
	Consider  the Dirichlet semigroup $P_t^D$ generated by $\L:=\DD+Z$ on $\bar D$, where $D$ is a regular open domain in a Riemannian manifold, $\DD$ is the Laplacian and $Z$ is a $C^1$ vector field on $\bar D$.
	We have
	$$P_t^Df(x)=\E[\1_{\{t<\tau(X^x)\}}f(X_t^x)],\ \ t\ge 0,\ x\in \bar D, \ f\in \B_b(\bar D),$$
	where $\B_b(\bar D)$ is the space of all bounded measurable functions on $\bar D$,  $X_t^x$ is the diffusion process generated by $\L$ with $X_0^x=x$, and
	$$\tau(X^x):=\inf\{t\ge 0: X_t^x\in\pp D\}$$
	is the first hitting time of $X_t^x$ to the boundary $\pp D$. Here and in the following,  we set $\inf \emptyset=\infty.$
	By using the local version of Bismut formula established in \cite{A},
	the following type gradient estimate is derived in \cite{AW}:
	\beq\label{G} |\nn P_t^Df(x)|\le \ff {c\|f\|_\infty}{\ss t\, (\rr_\pp(x)\land 1)},\ \ \ t\in (0,T],\ x\in D,\ f\in \B_b(\bar D),\end{equation}
{where $\rr_\pp$ is the Riemannian  distance function to the boundary $\pp D$}, and $c>0$ is a constant explicitly determined by local geometry quantities of $D$,
including the dimension,   curvature  bounds of the generator in $D$, and the second fundamental form on $\pp D$.  Below we formulate an equivalent inequality  \eqref{G'} to  \eqref{G}, where the    point $x\in D$ is replaced by    sub-probability measures on $D$,     
which is essential   for killed DDSDEs.

For any probability measure   $\mu$ on $\bar D$,  its restriction $\mu|_D$ to $D$ is a sub-probability measure on $D$, i.e.
$$\mu|_{D}\in \scr P_D:=\big\{\nu:\ \nu\ \text{is\ a\ measure\ on\ }D,\ \nu(D)\le 1\big\}.$$
Let  $P_t^{D*}\mu\in\scr P_D$ be defined as
$$(P_t^{D*}\mu)(A):= \mu(P_t^D \1_A)=\int_D P_t^D\1_A\d\mu= \int_D P_t^D\1_A\d\mu|_D,\ \ \ A \in  \B_D,\ t\in (0,T],$$
where $\B_D$ is the Borel $\si$-field on $D$.
Since   $P_t^{D*}\mu$ only depends on $\mu|_D$,
$P_t^*$ is well-defined on $\scr P_D$, i.e. for any  $\mu\in \scr P_D$,  $P_t^{D*}\mu$ is identified  as $P_t^{D*}\bar \mu $ for a probability measure $\bar\mu$ on $\bar D$ with $\bar\mu|_D=\mu$.
So,
\beq\label{PO}   (P_t^{D*}\mu)(f):= \int_D f \d (P_t^{D*}\mu) = \int_{D} P_t^Df \d\mu,\ \ \ t\in (0,T], \ \mu\in \scr P_D,\ f\in \B_b(\bar D).\end{equation}
Let
$$\rr(x,y):= \big\{|f(x)-f(y)|:\ f\in C^1(D),\ |(\rr_\pp\land 1)\nn f|\le 1\big\},$$
which is  the intrinsic distance in $D$  induced by the Riemannian metric with weight $(\rr_\pp\land 1)^{-2}$, see \cite{AW} for details. The associated truncated $1$-Wasserstein distance is defined as
$$\hat\W_1^\rr(\mu,\nu):= \inf_{\pi\in \C_D(\mu,\nu)} \int_{\bar D\times \bar D} [\rr(x,y)\land 1]\pi(\d x,\d y),$$
where $\pi\in \C_D(\mu,\nu)$ means that $\pi$ is a probability measure on $\bar D\times \bar D$ such that
$$\pi(A \times \bar D) =\mu(A),\ \ \pi(\bar D\times A)=\nu(A),\ \ \ \ A\in  \B_D.$$
Then the gradient estimate
\eqref{G} is equivalent to
$$|P_t^Df(x)-P_t^Df(y)|\le \ff {c\rr(x,y)} {\ss t} \|f\|_\infty,\ \ \ t\in (0,T],\ f\in \B_b(\bar D),\  x,y\in \bar D.$$
This together with $|P_tf(x)-P_tf(y)|\le 2\|f\|_\infty$ implies that for some constant $c'>0$,
$$|P_t^Df(x)-P_t^Df(y)|\le \ff {c'[1\land \rr(x,y)]} {\ss t} \|f\|_\infty,\ \ \ t\in (0,T],\ f\in \B_b(\bar D),\ x,y\in \bar D.$$
Hence, \eqref{G} implies
\beq\label{G'} \|P_t^{D*}\mu-P_t^{D*}\nu\|_{var} \le \ff{c'}{\ss t} \hat\W_1^\rr(\mu,\nu),\ \ t>0,\ \ \mu,\nu\in \scr P_D.\end{equation}
On the other hand, it is trivial that \eqref{G'} implies \eqref{G} for $c=c'.$ So, we may regard \eqref{G'} as gradient estimate on $P_t^D$. This is essential for killed DDSDE, since in this case  
the distributions for the solutions with deterministic initial values are   not enough to determine those with random initial values.

From now on, let $D$ be a connected, {possibly unbounded,} open domain in $\R^d$ with boundary $\pp D$, let $\bar D:=D\cup \pp D$. For any random variable $\xi$ on $\bar D$, its restricted distribution   to $D$ is defined as
$$\L_\xi^D := \P(\xi \in \cdot)\in \scr P_D.$$   
When different probability spaces are considered, we denote $\L^D_\xi$ by  $\L^D_{\xi|\P}$ to emphasize the underlying probability $\P$.

According to \cite{FMP}, for the solution of a regular killed SDE on $\bar D$,  the associated Dirichlet semigroup $P_t^D$ on $\B_b(\bar D)$ satisfies (see Lemma \ref{LM} below)
\beq\label{*G} \|\nn P_t^Df\|_\infty\le \ff{c}{\ss t}\|f\|_\infty,\ \ t\in (0,T],\ f\in \B_b(\bar D).\end{equation}
By letting $P_t^Df|_{D^c}=0$ for $t>0$, this is equivalent to
$$|P_t^D f(x)-P_t^Df(y)|\le \ff {c\|f\|_\infty}{\ss t}|x-y|,\ \ \ t\in (0,T],\ f\in \B_b(\bar D),\ x,y\in \bar D,$$
and hence implies that for some possibly larger constant $c'>0$
\beq\label{R}\|P_t^{D*}\mu-P_t^{D*}\nu\|_{var}\le \ff {c'}{\ss t} \hat\W_1(\mu,\nu),\ \ t\in (0,T], \ \mu,\nu\in \scr P_D,\end{equation}
where $\hat\W_1$ is defined as $\hat\W_1^\rr$ with $\rr(x,y)$ replaced by $|x-y|$:
$$\hat\W_1(\mu,\nu):= \inf_{\pi\in \C_D(\mu,\nu)} \int_{\bar D\times \bar D} [|x-y|\land 1]\pi(\d x,\d y).$$
On the other hand,  since $\rr_\pp\land 1\le 1$,
it is clear that
$\rr(x,y)\ge |x-y|$.
So, 
$ \hat \W_1 \le \hat\W_1^\rr,$ hence \eqref{R} is stronger than \eqref{G'}.

In this paper, we intend to establish the estimates  \eqref{G'} and \eqref{R} for  the following type killed DDSDE on $\bar D$, in the regular and singular cases respectively:
\begin{align}\label{E1}
\d X_t=\1_{\{t<\tau(X)\}}\left\{b_t(X_t,\L^D_{X_t})\d t+\si_t(X_t)\d W_t\right\},\ \ t\in[0,T],
\end{align}
where 	$T>0$ is a fixed time,  $W_t$ is the $m$-dimensional Brownian motion on a probability base (complete filtration probability space) $(\Omega,\{\F_t\}_{t\in [0,T]},\P)$,
$$
\tau(X):=  \inf\{t\in[0,T]: X_t\in\pp D\},
$$
and
$$b: [0,T]\times D\times \scr P_D\to \R^d,\ \ \si: [0,T]\times D\to \R^d\otimes\R^m$$
are measurable.

\beg{defn} \beg{enumerate} \item[$(1)$]
A continuous adapted process $(X_t)_{t\in [0,T]}$ on $\bar D$ is called a solution of \eqref{E1}, if $\P$-a.s.
$$\int_0^{T\land \tau(X)} \big\{|b_t(X_t,\L^D_{X_t})|+\|\si_t(X_t)\|^2\big\}\d t<\infty,$$ and
$$X_t=X_0+\int_0^{t\land \tau(X)} \big\{b_s(X_s, \L^D_{X_s})\d s+ \si_s(X_s)\d W_s\big\},\ \ t\in [0,T].$$
We call \eqref{E1} strongly  well-posed, if for any $\F_0$-measurable random variable $X_0$ on $\bar D$, the SDE has a unique solution.
\item[$(2)$]
A couple $(X_t, W_t)_{t\in [0,T]}$ is called a weak solution of \eqref{E1}, if there exists a complete filtration probability space $(\Omega,\{\F_t\}_{t\in [0,T]},\P)$ such that
$W_t$ is the $m$-dimensional Brownian motion and $X_t$ solves \eqref{E1}. We call \eqref{E1} weakly  well-posed, if for any $\mu\in \scr P_D$, there exists a weak solution $(X_t,W_t)_{t\in [0,T]}$ of $\eqref{E1}$
with respect to a probability base $(\Omega,\{\F_t\}_{t\in [0,T]},\P)$ such that  $\L_{X_0|\P}^D=\mu$, and if there is another weak solution $(\tt X_t,\tt W_t)_{t\in [0,T]} $
with respect to   probability base $(\tt\Omega,\{\tt\F_t\}_{t\in [0,T]},\tt\P)$ such that  $\L_{X_0|\tt\P}^D=\mu$, then $\L_{X_t|\P}^D=  \L_{\tt X_t|\tt\P}^D $ for $t\in [0,T].$
\item[$(3)$]
We call \eqref{E1} well-posed, if it is both strongly  and weakly  well-posed.  \end{enumerate} \end{defn}

When \eqref{E1} is well-posed,  we let $P_t^{D*}\mu=\L_{X_t}^D$ for the solution with $\L_{X_0}^D=\mu$, and study the
the regularity of  the maps
$$\scr P_D\ni \mu \mapsto P_t^{D*}\mu\in \scr P_D,\ \  t\in (0,T]$$
by establishing the estimate \eqref{G'} and \eqref{R} for $P_t^{D*} $ in place of $P_t^*$.

To measure the   singularity, we recall locally integrable functional spaces introduced in \cite{XXZZ}.
For any $t>s\ge 0$ and $p,q\in (1,\infty)$, we write $f\in  \tilde{L}_p^q([s,t])$  if $f: [s,t]\times\R^d\to \R$ is measurable with
$$\|f\|_{\tilde{L}_p^q([s,t])}:= \sup_{z\in\R^d}\bigg\{\int_{s}^t \bigg(\int_{B(z,1)}|f(u,x)|^p\d x\bigg)^{\ff q p} \d u\bigg\}^{\ff 1q}<\infty,$$ where $B(z,1):=\{x\in \R^d: |x-z|\le 1\}$ is the unit ball centered at point $z$.
When $s=0$, we simply denote
\begin{equation*} \tt L_p^q(t)=\tt L_p^q([0,t]),\ \ \|f\|_{\tilde{L}_p^q(t)}=\|f\|_{\tilde{L}_p^q([0,t])}.\end{equation*}
We will take $(p,q)$ from the   space
\begin{equation*} \scr K:=\Big\{(p,q) \in (2,\infty)\times (2,\infty):\   \ff{d}p+\ff 2 q<1\Big\}.\end{equation*} 
In Section  2, we study   the killed SDE where $b_t(x,\mu)=b_t(x)$ does not depend on the distribution variable $\mu$,  by first establishing Bismut formula   in the regular case, then deriving gradient estimate in the singular case. As applications, 
{in Section 3, the estimate \eqref{G'} is presented for singular killed DDSDEs, and the stronger estimate \eqref{R}  for regular killed DDSDEs.}
It is worthy indicating that the Bismut formula derived in Section 2 is new,  since the existing formula in \cite{A,AW} only works for time homogeneous models with diffusion coefficient at least $C^2$-smooth, so that
the required curvature condition is well-defined.

\section{Bismut formula and gradient estimate }

In this section,  we  let $b_t(x,\mu)=b_t(x)$ be independent of $\mu$, and consider the classical SDE on $\R^d$:
\begin{equation}\label{E1'} \d X_t= b_t(X_t)\,\d t+\si_t(X_t)\,\d W_t,\ \ t\in [0,T].\end{equation}
Under  the following assumption, the well-posedness and Bismut formula have been derived in \cite{HRW} for this SDE, see \cite{W23b, XXZZ} for earlier results where $k_1=1$.

\emph{\beg{enumerate} \item[{\bf (A)}]  Let $a:=\si\si^*$,  $b=b^{(0)}_t+\sum_{i= 1 }^{k_1} b_t^{(i)}$ for some $k_1\in \mathbb N$ and  measurable maps
$$b^{(i)}: [0,T]\times\R^d\to \R^d,\ \  0 \le i\le k_1, $$  such that the following  conditions hold.
\item[$(1)$]  $ \|a_t(x)\|+\| a_t(x)^{-1}\|+ \|\nn b_t^{(0)}(x)\|+|b_t^{(0)} (0)|$ is  uniformly  bounded in  $(t,x)\in[0,T]\times\R^d$, and
$$\lim_{\vv\downarrow 0} \sup_{t\in [0,T], |x-y|\le \vv} \|a_t(x)-a_t(y)\| =0. $$
\item[$(2)$] There exist $k_2 \in \mathbb N$, $\{(p_i,q_i), (p_j', q_j')\}_{1\le i\le k_1,\ 1\le j\le k_2}\subset \scr K$, and $\big\{0\le f_j \in \tt L_{q_j'}^{p_j'}  (T)  \big\}_{1\le j\le k_2}$ such that
$$\sup_{1\le i\le k_1}  \|b^{(i)}\|_{\tt L_{q_i}^{p_i} (T) }<\infty,\ \ \ \|\nn a\|\le \sum_{j=1}^{k_2} f_j.$$
\end{enumerate}}

According to   \cite[Proposition 5.2]{HRW}, which extends \cite[Theorem 1.1.3]{RW24} with $k_1=1$,   {\bf (A)} implies that  the SDE  \eqref{E1'} is well-posed, i.e. for any $x\in \R^d$, it has a unique solution  $X_t^x$ with $X^x_0=x$. Moreover,   for any $x,v\in \R^d$, the derivative process
$$\nn_vX_t^x:=\lim_{\vv\downarrow  0} \ff{X_t^{x+\vv v}-X_t^x}\vv,\ \ t\in [0,T]$$ exists in $L^k(\OO\to C([0,T];\R^d),\P)$ for any $k\in [1,\infty)$, and
there exists a constant $c(k)>0$ such that
\begin{equation}\label{DFF0} \sup_{x\in \R^d} \E\Big[\sup_{t\in [0,T]}|\nn_v X_t^x|^k\Big]\le c(k) |v|^k,\ \ \ x,v\in\R^d.\end{equation}
Now, let
$$
\tau(x)=\tau(X^x):=  \inf\big\{t\in [0,T]:\ {X_t^x}\in\pp D\big\},\ \  x\in \bar D.$$
Then the associated Dirichlet semigroup is defined as
$$P_t^D f(x):= {\E}\big[\1_{\{t<\tau(x)\}} f(X_t^x)\big],\ \ \ t\in [0,T],\ x\in \bar D,\ f\in \B_b(\bar D).$$
In general, for any $(s,x)\in [0,T)\times \R^d$, let $(X_{s,t}^x)_{t\in [s,T]}$ be the unique solution to \eqref{E1'} from time $s$ with $X_{s,s}^x=x$,  and define
\beq\label{TS} \beg{split} &\tau(s,x):= \inf\big\{t\in [s,T]:\ X_{s,t}^x\in\pp D\big\},\\
& P_{s,t}^D f(x):= {\E}\big[\1_{\{t<\tau (s,x)\}} f(X_{s,t}^x)\big],\ \ \ t\in [s,T],\ x\in \bar D,\ f\in \B_b(\bar D).\end{split} \end{equation}
Then $P_t^D =P_{0,t}^D$ for $t\in [0,T],$ and   the Markov property of the solution  implies  the flow property
\beq\label{FL} P_{s,t}^D= P_{s,r}^D P_{r,t}^D,\ \ \ 0\le s\le r\le t\le T.  \end{equation}
In the following, we first establish a local version Bismut formula for $P_t^D$ in the regular case,
then derive gradient estimate on $P_t^D$ in the singular case.

\subsection{Bismut formula in the regular case}

We first recall some notations.
Let $M$ be a Riemannian manifold possibly with a boundary $\pp M$. For any $\aa\in (0,1)$ and $k\in \Z_+$, let $ C^{k+\aa}_{loc} (M)$ be the set of all real functions $g$ on $M$ such that $\nn^l g$ exists for $0\le l\le k$, where $\nn^0g:=g,$ and $\nn^{k} g$ is locally $\aa$-H\"older continuous.
We write $g\in C_b^{k+\aa}(M)$, if $\nn^lg$ is bounded for any $0\le l\le k$, and $\nn^k g$ is globally $\aa$-H\"older continuous on $M$.

Let $I\subset [0,T]$ and $M\subset \R^d$.
For any $\aa_1,\aa_2\in (0,1]$, we write $f\in C_{b}^{\aa_1,\aa_2}(I\times M),$ if $f\in C_b(I\times M)$  such that for  a constant $c>0$
$$|f_t(x)-f_s(y)|\le c \big(|t-s|^{\aa_1}+|x-y|^{\aa_2}\big),\ \ t,s\in I,\ x,y\in M.$$
We write $f\in C_{loc}^{\aa_1,\aa_2}(I\times M)$ if  $f\in  C_b^{\aa_1,\aa_2}(I'\times M')$ holds for any compact sets $I'\times M'\subset I\times M$.
Moreover, for any $k_1,k_2\in \Z_+,$ we write $f\in C_{loc}^{k_1+\aa_1,k_2+\aa_2}(I\times M)$ if   { $\pp_t^{l_1}f, \nn^{l_2} f\in C_{loc}^{\aa_1,\aa_2}(I\times M)$}
holds  for any integers $0\le l_1\le k_1$ and $0\le l_2\le k_2$, and $f\in C_b^{k_1+\aa_1,k_2+\aa_2} (I\times M)$ if    {  $\pp_t^{k_1}f, \nn^{k_2} f\in C_b^{\aa_1,\aa_2} (I\times M)$ } holds  for any integers $0\le l_1\le k_1$ and $0\le l_2\le k_2$.

For any $k\in \Z_+$ and $\aa\in (0,1)$, we say that $\pp D$ is uniformly of class $C^{k+\aa}$, if it has a unique inward unit normal vector field $N$ such that for some constant $r_0>0$, the polar coordinates
$$\pp D\times [0,r_0]\ni (z,r)\mapsto z+rN(z)\in D_{r_0}:= \big\{x\in\bar D:\ \rr_\pp(x):={\rm dist}(x,\pp D) \le r_0\big\}$$
is bijective, and each component belongs to $C_b^{k+\aa}(\pp D\times [0,r_0]).$ In this case, we have $\rr_\pp\in C_b^{k+\aa}(D_{r_0}).$

To establish a local version Bismut formula for $P_t^D$, we make the following assumption.

\emph{\begin{enumerate}
\item [{\bf (B)}] Let $a:=\si\si^*$.  
There exist   constants $\alpha \in(0,1)$ and $K>0$  such that
\item[$(1)$]  $\pp D$ is uniformly of class $C^{3+\aa}$,
or uniformly of class $C^{2+\aa}$ and $(a_t,b_t)=(a,b)$ does not depend on $t$.
\item[$(2)$] { Either  $(a_t,b_t)=(a,b)$ do not depend on $t$, or  
\beq\label{F2}   \|1_{D_{r_0} } b\|_\infty\le K\ \text{\ for \ some \ }r_0>0.\end{equation}}
\item[$(3)$] All components of $b$ and $a$ belong to $C_{loc}^{\aa/2, 1} ([0,T]\times \R^d)$ such that
\beq\label{F1} \beg{split} & \|a\|_\infty +\|a^{-1}\|_\infty+  \|\nn \si\|_\infty^2   \le K,\\
&\<\nn_v b_t(x),v\>\le K|v|^2,\ \ \ (t,x,v)\in [0,T]\times\R^d\times\R^d.\end{split}\end{equation}
\end{enumerate}}

We have the following result.

\beg{thm}\label{T1} Assume {\bf(B)}. For $(t,x)\in (0,T]\times D$,    let $(\bb_s)_{s\in [0,t]}$ be an adapted absolutely continuous real process such that
\beq\label{BB1} \bb_0=1, \ \ \bb_s=0\ \text{for\ } s\geq \tau(x)\wedge t,\end{equation}  and     for some $\vv>0$
\beq\label{BB2} C_{\bb,\vv,{t}}:=\E\left(\int_{0}^{t\wedge\tau(x)}|\bb_s'|^2\,\d s \right)^{\ff{1+\vv} 2} <\infty,\end{equation}
then for any 
$v\in \R^d$ and $f\in \B_b(\bar D)$,
\begin{equation}\label{BSM1} \nn_v P_t^Df(x)= {-}\E\bigg[\1_{\{t<\tau(x)\}}f(X_t^x) \int_0^{t} \bb_s'\Big\<(\si^*a^{-1})_s(X_s^{x}) \nn_vX_s^x,  \d W_s\Big\> \bigg].\end{equation}
Consequently, for any $p\in  (1+\vv^{-1},\infty]$,
there exists a constant $c(p,\vv)>0$ such that for any $t\in(0,T],\,x\in D$ and $f\in \B_b(\bar D)$,
\begin{equation}\label{GRD}    \left|\nn P_t^Df(x)\right|:= \limsup_{y\to x} \ff{|P_t^Df(y)- P_t^Df(x)|}{|x-y|}\le   c(p,\vv) C_{\bb,\vv,{t}}^{\ff 1{1+\vv}}\ (P_t^D|f|^p(x))^{\ff 1 p}. \end{equation} 
\end{thm}

{ \beg{rem} We would like to make comments on the local Bismut formula  \eqref{BSM1} and the condition \eqref{F2}.
\beg{enumerate} 
\item[$(1)$] When $(a_t,b_t)=(a,b)$ does not depend on $t$ and is $C^2$-smooth, the local Bismut formula of type \eqref{BSM1}  was first established and applied   in \cite{A,AW}.
The above result extends the main results in these papers to   time dependent and less regular coefficients, which will be further applied to the singular case for coefficients satisfying assumption {\bf (A)}, see 
Theorem \ref{T2} below.  
\item[$(2)$] The condition \eqref{F2}  is used to ensure the gradient estimate \eqref{U2} below due to \cite[Theorem 5.3]{AL}. According to \cite[Theorem 1.3]{FMP}, this estimate holds when $(a_t,b_t)=(a,b)$ does not depend on $t$,
since \eqref{F1} implies   the Lyapunov condition
\beq\label{LYY} \ff 1 2{\rm tr}\big\{ a\nn^2 V\big\}+b\cdot\nn V\le \ll_0 V\end{equation} 
on $\bar D$ for some $\ll_0>0$ and $V(x)=|x|^2$.
\end{enumerate} 
\end{rem} }

To prove this theorem, we first present the following lemma.

\beg{lem}  \label{LM} Assume {\bf (B)}.
\beg{enumerate}
\item[$(1)$] There exists a constant $c>0$ depending only on $d, T,   r_0,  K$ and $\pp D$ such that
\beq\label{U2} \|\nn P_{s,t}^D f\|_\infty\le \ff {c \|f\|_\infty} {\ss{t-s}} ,\ \ \  0\le s<t\le T,\ f\in \B_b(\bar D).\end{equation}
\item[$(2)$] For any $f\in \B_b(\bar D)$ and $t\in (0,T],$ we have  $P_{\cdot,t}^D f\in C_{loc}^{1+\ff\aa 2, 2+\aa}([0,t)\times \bar D),$ and
\beq\label{U1} \pp_s P_{s,t}^Df= -L_s P_{s,t}^Df,\ \ \ s\in [0,t),\end{equation}
where $$L_{s}:= {\rm tr}\{a_s\nn^2\}+ b_s\cdot\nn,\ \ s\in [0,t).$$ 
\end{enumerate}
\end{lem}
\beg{proof} (1)  Let $P_t^Df|_{D^c}=0$ for $t>0$. The gradient estimate \eqref{U2} is equivalent to
$$ |P_{s,t}^D f(x)-P_{s,t}^Df(y)|  \le \ff {c |x-y|} {\ss{t-s}}\|f\|_\infty,\ \ \  0\le s<t\le T,\,f\in\B_b(\bar D),\, x,y\in \bar D.$$
So, by the monotone class theorem, we  only need to prove  \eqref{U2} for $f\in C_b(\bar D)$.

If,   in addition to {\bf (B)},  all components of  $a$ and $ b$ belong to $ C_{loc}^{\aa/2, 1+\aa}([0,T]\times \R^d)$, then \eqref{U2} is known. Indeed,
when $\pp D$ is uniformly of class $C^{3+\aa}$,   \eqref{U2} is derived in \cite[Theorem 5.3]{AL} for $f\in C_b(\bar D)$;  when $b_t$ and $a_t$ do not depend on $t$  and $\pp D$ is uniformly of class $C^{2+\aa}$,
\eqref{U2} can be found in
\cite[Theorem 1.3]{FMP}, {where the constant $c$ only   depends on $K, T$ and $\pp D$}, see also   \cite{L16}.

In general,   \eqref{U2} can be derived by using mollifier approximations as follows.   Let $0\leq g\in C^\8_0(\R^{d+1})$ with support contained in $\{(r,x)\in \R\times \R^{d}: |(r,x)|\leq1\}$ such that $\int_{\R^{d+1}}g(r,x)\,\d r\d x=1$. For any $n\geq1$, let $g_n(r,x)=n^{d+1}g(nr,nx)$  and  define
\begin{align*}
&a^n_t(x):= \si^n_t(x)\si^n_t(x)^*,\ \ \si^{n}_t(x):=\int_{\R\times\R^d}\si_{(t+s)^+\wedge T}(x+y)g_n(s,y)\,\d s\,\d y,\\
&b^{n}_t(x):=  \int_{\R\times\R^d} b_{(t+s)^+\wedge T}(x+y)g_n(s,y)\,\d s\,\d y,\  \  (t,x)\in \R\times\R^d.
\end{align*}
Then $b^n,a^n\in C^\infty([0,T]\times \R^d)$, \eqref{F1} holds for $(b^n, a^n)$ in place of $(b,a)$ with the same constant $K$,  and when $n\ge n_0$ for a  larger enough constant $n_0\ge 1$,
\eqref{F2} holds for $( b^n, r_0/2)$    in place of $(b,r_0)$ with the same constant $K$.
So, there exists a constant $c>0$  depending only on $d, T,   r_0,  K$ such that
$$ \|\nn P_{s,t}^{D,n} f\|_\infty\le \ff {c \|f\|_\infty} {\ss{t-s}} ,\ \ \  0\le s<t\le T,\ f\in \B_b(\bar D),\  n\ge n_0,$$
where $P_{s,t}^{D,n}$ is defined as in \eqref{TS} for the solution to $X^{n,x}$ to the following SDE in place of  $X^x$:
\begin{align*}
\d X^{n,x}_{s,t}=b^n_t(X^{n,x}_{s,t})\d t+\si_t^n(X^{n,x}_{s,t})\d W(t),\ \ 0\leq s\leq t\leq T, \ \ X_{s,s}^{n,x}=x\in D.	\end{align*}
Noting that
$$\lim_{n\to\infty} P_{s,t}^{D,n}f = P_{s,t}^D f,\ \  0\le s<t\le T,\, f\in C_b(\bar D),$$
see the proof of \eqref{LMT} in a more general setting where $b$ may be singular, we derive \eqref{U2}.

(2)  For any fixed $t\in(0,T]$ and $s\in[0,t)$,   since     $P_{r,t}^Df\in C_b(\bar D)$ for  $f\in\B_b(\bar D)$  and $r\in (s,t)$ due to \eqref{U2}, by the flow property
\eqref{FL} we only need to prove \eqref{U1}  for  $f\in C_b(\bar D).$
In this case,
by \cite[Theorem 3.4]{AL},  the following PDE has a unique bounded classical solution $u\in C^{1+\aa/2, 2+\aa}_{loc} ((0, t]\times \bar D)$:
\beq\label{*E}  \beg{split} &\pp_r u_r(x)= L_{t-r}u_r(x),\ \ \ (r,x)\in (0,t]\times D,\\
&u_0=f,\ u_r|_{\pp D}=0,\ \ r\in (0,t].\end{split}\end{equation}
By It\^o's formula,
$${\d}u_{t-r}(X_{s,r}^x)= \big\<\nn u_{t-r}(X_{s,r}^x), \si_r(X_{s,r}^x)\d W_r\big\>,\ \  x\in\bar D,\, r\in [s,t)\cap [s,  \tau(s,x)].$$
Combining this and the fact that   for any $(s,x)\in [0,T)\times\bar D$,
$$u_{t-t\land  \tau(s,x)} (X^x_{s,t\land \tau(s,x)})= u_{t-\tau(s,x)}(X^x_{s,\tau(s,x)})=0\ \text{if}\ t \geq \tau(s,x),$$  we derive
\beg{align*}&u_{t-s}(x)= \E [u_{t-s}(X_{s,s}^x)]= \E[u_{t-t\land \tau(s,x)}(X_{s,t\land \tau(s,x)}^x) ]\\
&= \E[1_{\{t < \tau(s,x)\}} f(X_{s,t}^x)] = P_{s,t}^D f(x),\ \  x\in\bar D,\ 0\le s\le t.  \end{align*}
Then \eqref{*E} implies \eqref{U1}.
\end{proof}

\begin{proof}[Proof of Theorem \ref{T1}]    For any fixed $x\in D$, we simply denote $X_t^x=X_t$ and $\tau(x)=\tau$.  Let $f\in \B_b(\bar D).$

(a) Below, we prove that \eqref{GRD} can be deduced from \eqref{BSM1}. By \eqref{F1},  for any $v\in\R^d$, the derivative process $v_t:= \nn_v X_t$ solves the SDE
$$\d v_t= (\nn_{v_t} b_t)(X_t)\d t + (\nn_{v_t}\si)(X_t)\d W_t,\ \ v_0=v,\ t\in [0,T],$$
and
$$\d |v_t|^2  =\big\{ 2\<(\nn_{v_t} b_t)(X_t), v_t\>+ \|\nn_{v_t}\si\|^2(X_t)\big\}\d t + \d M_t
\le 3K|v_t|^2\d t+  \d M_t,\ \ t\in [0,T],$$
where $\d M_t:= 2\big\<v_t, (\nn_{v_t}\si)(X_t)\d W_t\big\>$ satisfies
$$\d \<M\>_t\le 4K |v_t|^2\d t.$$
By the Burkholder-Davis-Gundy  inequality, for any $k\in [1,\infty)$, we find a constant $c(k)>0$ such that \eqref{DFF0}  holds.
Combining this with  H\"older's inequality,    for any  $t\in[0,T]$,  $p> 1+\vv^{-1}$ and $q= \ff{(p-1)(1+\vv)}{\vv p-1-\vv},$     we obtain
\beq\label{O1}\beg{aligned}& A_t:= \E\bigg[\bigg(\int_0^{t\wedge\tau} \big|\bb_s' (\si^*a^{-1})_s(X_s^{x}) \nn_vX_s\big|^2 \d s\bigg)^{\ff p{2(p-1)}}\bigg] \\
&\le  \|\si^* a^{-1}\|_\infty^{\ff p{p-1}} \E\bigg[\bigg(\Big(\sup_{s\in [0,t]} |\nn_v X_{s}|^{\ff p{p-1}} \Big) \bigg(\int_0^{t\wedge\tau} \big|\bb_s' \big|^2 \d s\bigg)^{\ff p{2(p-1)}}\bigg]\\
&\le  \|\si^* a^{-1}\|_\infty^{\ff p{p-1}}  \bigg(\E\bigg[ \sup_{s\in [0,t]} |\nn_v X_s|^{\ff {pq}{p-1}} \bigg]\bigg)^{\ff{1}{q}}
\bigg[\E\bigg(\int_0^{t\land\tau} |\bb_r'|^2 \d r\bigg)^{\ff{1+\vv} 2}\bigg]^{\ff {p}{(p-1)(1+\vv)}}<\infty\end{aligned}\end{equation} due to   \eqref{DFF0},   \eqref{F1}  and \eqref{BB2}.
By \eqref{BSM1},    H\"older's inequality and the Burkholder-Davis-Gundy  inequality, 
we find  a constant $c_p>0$ such that  for any  $t\in(0,T]$, 
\begin{equation}\label{GRD0}
\begin{aligned} & |\nn_v P_t f(x)|= \bigg|\E\bigg[\1_{\{t<\tau\}}f(X_t) \int_0^{t\wedge\tau} \bb_s'\Big\<(\si^*a^{-1})_s(X_s) \nn_vX_s,  \d W_s\Big\> \bigg]\bigg| \\
& \le c_p  \big(P_t^D |f|^p(x)\big)^{\ff 1 p}  A_t^{\ff{p-1}p} \\
&\le c_p\|\si^* a^{-1}\|_\infty \big(P_t^D |f|^p(x)\big)^{\ff 1 p} \E\bigg( \sup_{s\in [0,t]} |\nn_v X_s|^{\ff {pq}{p-1}} \bigg)^{\ff{p-1}{pq}}
\bigg[\E\bigg(\int_0^{t\land\tau} |\bb_r'|^2 \d r\bigg)^{\ff{1+\vv}  2 }\bigg]^{\ff 1{1+\vv}}.\end{aligned}
\end{equation}
So, \eqref{GRD} follows from  \eqref{DFF0} and \eqref{F1}. 

(b)   We now  establish \eqref{BSM1}.  For any fixed $t\in(0,T]$,  by  \eqref{BB1} and  \eqref{O1},
$$M_s:=\int_0^s \bb_r'\big\<(\si^* a^{-1})_r(X_r)\nn_{v} X_r,\d W_r\big\>= \int_0^{s\land\tau} \bb_r'\big\<(\si^* a^{-1})_r(X_r)\nn_{v} X_r,\d W_r\big\>,\ \ s\in  [0,  t] $$
is a uniformly integrable martingale.     By  \eqref{U2}, \eqref{U1}  and It\^o's formula,
\beq\label{U3}   P_{s,t}^Df(X_s)-P_t^Df(x)  =\int_0^s \big\<\nn P_{r,t}^D f(X_r), \si_r(X_r)\d W_r\big\>=: N_s,\ \ s\in [0,t)\cap[0,  \tau] \end{equation}  is a martingale.
So,
$$I(s) := P_{s,t}^D f(X_s) \int_0^s \bb_r'\big\<(\si^* a^{-1})_r(X_r)\nn_{v} X_r,\d W_r\big\>,\ \ s\in [0,t)\cap[0,  \tau] $$
is a uniformly integrable process with
\beq\label{U4}  I(s)=(N_s+ P_t^Df(x)) M_s  = I_1(s)+I_2(s)+I_3(s), \ \ \ s\in [0,t){\cap} [0,\tau], \end{equation}
where
\beg{align*}& I_1(s):= \int_0^s \bb_r'\big\<\nn_v X_r,\nn P_{r,t}^D f(X_r)\big\>\d r,
\ \ s\in [0,t),
\\
&I_2(s):=   \int_0^s
P_{r,t}^D f(X_r) 
\bb_r' \big\<(\si^* a^{-1})_r(X_r)\nn_{v} X_r,\d W_r\big\>,
\ \ s\in [0,t),
\\
& I_3(s):=   \int_0^s M_r\d N_r,
\ \ \ \ s\in [0,t)\cap[0, \tau].
\end{align*}
By \eqref{DFF0}, \eqref{F1},  \eqref{U1} and the condition on $\bb_s$, we see that $I_2(s)$ and $I_3(s)$ are   martingales for $s\in  [0,t)\cap[0,  \tau].$
Moreover, by $\bb_0=1$, the chain rule,   integration by parts formula, and \eqref{U3}, we obtain
\beq\label{I1} \beg{split} I_1(s)&=  \int_0^s \bb_r'\nn_v \big[P_{r,t}^D f(X_r)\big]\d r  \\
&= \bb_s\nn_v\big[P_{s,t}^D f(X_s)\big]- \nn_v P_t^Df(x) -\int_0^s\bb_r \nn_v\big[\d P_{r,t}^D f(X_r)\big]\\
&= \bb_s\nn_v\big[P_{s,t}^D f(X_s)\big]- \nn_v P_t^Df(x) - \tt M_s,\ \ s\in [0,t),\end{split}
\end{equation}
where, by Lemma \ref{LM}  and the boundness of $\bb_r$ and $\nn \si$,
$$\tt M_s:=  \int_0^s\bb_r \Big[\Hess_{P_{r,t}^D f(X_r)}(\nn_v X_r, \si_r(X_r) \d W_r) +\big\<\nn P_{r,t}^D f(X_r),\nn_{\nn_vX_r}(\si_r)(X_r) \d W_r\big\> \Big],\ \ s\in [0,t)$$
is a local martingale, i.e. for
$${\hat\tau_n}:=\inf\{t\ge 0:\ |X_t|\ge n\},\ \ n\ge 1,$$
we have ${\hat\tau_n}\to\infty$ as $n\to\infty$,  and for each $n\ge 1$,
$$[0,t) \ni s\to \tt M_{s\land{\hat\tau_n}} $$ is a martingale.
Combining \eqref{U4} with \eqref{I1}, and noting that $I_2(s)+I_3(s)$ is a martingale for
$s\in  [0,t)\cap[0,  \tau]$, we conclude that for any $n\ge 1$ and $s<t$,
\beg{align*} & 0 = I(0) = \E[I(s\land\tau\land{\hat\tau_n})-I_1(s\land\tau\land{\hat\tau_n})]\\
&= \E\Big(I(s\land\tau\land{\hat\tau_n}) +\nn_v P_t^Df(x)-
\bb_{s\land \tau\land{\hat\tau_n}} \nn_v \big[P_{s\land\tau\land{\hat\tau_n}, t} f(X_{s\land\tau\land{\hat\tau_n}})\big]\Big).\end{align*}
Combining this with
the uniform integrability of $\{I(s\wedge\tau)\}_{s\in[0,t)}$, the boundness of $\bb$,  and noting that \eqref{U2} implies
$$\sup_{  r\in [0,s]} \|\nn P_{r,t}^D f\|_\infty <\infty,\ \ s\in [0,t),$$
it follows from dominated convergence theorem that 
\beq\label{I2} \beg{split}\nn_v P_t^Df(x)&= \lim_{s\uparrow t, n\uparrow\infty} \bigg[\E\Big(-I(s\land\tau\land{\hat\tau_n})
+\bb_{s\land \tau\land{\hat\tau_n}} \nn_v \big[P_{s\land\tau\land{\hat\tau_n}, t} f(X_{s\land\tau\land{\hat\tau_n}})\big]\Big)\bigg]\\
&= \lim_{s\uparrow t } \bigg[\E\Big(-I(s\land\tau)
+\bb_{s\land \tau} \nn_v \big[P_{s\land\tau, t} f(X_{s\land\tau})\big]\Big)\bigg].\end{split}\end{equation}
Noting that $\P(\tau=t)=0 $ and $P_{\tau,t}^D f(X_\tau)=0$ for $\tau <t$, we have {$\P$-a.e.,}
$$P_{t\land\tau,t}^Df(X_{t\land\tau}) = 1_{\{t<\tau\}} f(X_t) + 1_{\{t>\tau\}} P_{\tau,t}^D f(X_\tau)=1_{\{t<\tau\}} f(X_t),$$
so that, by the uniform integrability of  $\{I(s\wedge\tau)\}_{s\in[0,t)}$ 
and the dominated convergence theorem,
\beq\label{I3} \lim_{s\uparrow t }  \E [I(s\land\tau)]= \E[I(t\land\tau)]
= \E \bigg[1_{\{t<\tau\}} f(X_t) \int_0^t \bb_r'\big\<\zeta_r(X_r)\nn_{v} X_r,\d W_r\big\>\bigg].\end{equation}
By  \eqref{DFF0}, \eqref{U2} and $\|f\|_\infty<\infty$, we find  constants $c_0,c_1>0$ independent of   $t$ and $s\in[0,t)$, such that 
\begin{equation}\label{I4}
\begin{aligned} &\E \big|\bb_{s\land \tau} \nn_v P^D_{s\land\tau, t} f(X_{s\land\tau})\big|
= \E \big|\bb_{s\land \tau}   \big\<\nn P^D_{s\land\tau, t} f(X_{s\land\tau}),\nn_v X_{s\land \tau}\big\>\big|\\
&\le \ff {c_0}{\ss{t-s}}   \Big(\E\big|\bb_{s\land \tau}\big|^{1+\vv}\Big)^{\ff 1 {1+\vv}} \left(\E | \nn_v X_{s\wedge\tau} |^{\ff{1+\vv}\vv} \right)^{\ff\vv{1+\vv}}\le \ff {c_1 |v| }{\ss{t-s}}   \Big(\E\big|\bb_{s\land \tau}\big|^{1+\vv}\Big)^{\ff 1 {1+\vv}}.\end{aligned}
\end{equation}
Moreover, by $\bb_{t\land\tau}=0 $  and H\"older's inequality, we have
$$|\bb_{s\land\tau}|^{1+\vv}= \bigg|\int_{s\land\tau}^{t\land\tau}\bb_r'\,\d {r}\bigg|^{1+\vv}\le (t-s)^{\ff{1+\vv}2} \bigg(\int_{s\land\tau}^{t\land\tau}|\bb_r'|^2\,\d {r}\bigg)^{\ff  {1+\vv} 2},\ \  s\in[0,t]$$
so that  \eqref{BB2}  implies
\beg{align*} \ff {1}{\ss{t-s}}   \Big(\E\big|\bb_{s\land \tau}\big|^{1+\vv}\Big)^{\ff 1 {1+\vv}}
\le      \bigg(\E\bigg[ \bigg(\int_{s\land\tau}^{t\land\tau} {|\bb_r'|^2\,\d r}\bigg)^{\ff{1+\vv}2}\bigg]\bigg)^{\ff 1 {1+\vv}}
\to 0\ \text{as}\ s\uparrow t.\end{align*}
This together with \eqref{I2}, \eqref{I3}  and \eqref{I4}   implies \eqref{BSM1} for $f\in \B_b(\bar D). $
\end{proof}

\subsection{Gradient estimate in the singular case}

\emph{\beg{enumerate} \item[{\bf (C)}]   Assume that {\bf (A)} and {\bf (B)}$(1)$-$(2)$ hold for some $r_0>0$.   \end{enumerate}}


\beg{thm}\label{T2} Assume {\bf (C)}, and let $P_t^D$ be associated with the SDE $\eqref{E1'}$. Then for any $p>2$ there exists a constant $c(p)>0$ such that for any $ t\in (0,T], \ x\in D$ and $ f\in \B_b(\bar D),$
\beq\label{GF}  |\nn P_t^Df(x)|
\le \ff {c(p)}  {\ss t\,(\rr_\pp(x)\land 1)}  \big(P_t^D|f|^p(x)\big)^{\ff 1 p}.\end{equation}
\end{thm}

\beg{proof} By induction in $k_1\in \mathbb N$ as in the proof of \cite[Proposition 5.2]{HRW},  we only need to prove for $k_1=1$.

(a)   To apply \eqref{GRD}, we first make regular approximation of the SDE. Let $0\leq g\in C^\8_0(\R^{d+1})$ with support contained in $\{(r,x)\in \R\times \R^{d}: |(r,x)|\leq1\}$ such that $\int_{\R^{d+1}}g(r,x)\,\d r\d x=1$. For any $n\geq1$, let $g_n(r,x)=n^{d+1}g(nr,nx)$  and
\begin{align*}
&b^{n}_t(x)=(b^{(0)}+  b^{1,n})_t(x),\ \ b_t^{1,n}(x):= \int_{\R\times\R^d} b^{(1)}_{(t+s)^+\wedge T}(x+y)g_n(s,y)\,\d s\,\d y,\\
& \si^{n}_t(x)=\int_{\R\times\R^d}\si_{(t+s)^+\wedge T}(x+y)g_n(s,y)\,\d s\,\d y,\ \
a^n_t(x):= \si^n_t(x)\si^n_t(x)^*,\ \  (t,x)\in \R\times\R^d.
\end{align*}
For any $n\geq1$ and any  $x\in D$,  denote $\tau_n(x)=\tau(X^{n,x})$, and $$P_t^{D,n}f(x):=\E\left[\1_{\{t<\tau(X^{n,x})\}}f(X_t^{n,x})\right],\ \ t\in[0,T],\,f\in\B_b(\bar D),$$
where $X_t^{n,x}$ solves
\begin{align*}
\d X^{n,x}_t=b^n_t(X^{n,x}_t)\d t+\si_t^n(X^{n,x}_t)\d {W_t},\ \ t\in[0,T], \ \ X_0^{n,x}=x\in D.	\end{align*}
Under {\bf (C)},    for any $n\geq1$,  {\bf(B)}  holds for $(b^n,a^n)$ in place of $(a,b)$. By
\eqref{GRD0} with $\vv=1$ for  $P_{t}^{D,n}$ in place of $P_t^D$,   for any $p>2$, we have
\beq\label{GRDn}  \beg{split}  \left|\nn_v P_t^{D,n}f(x)\right|\le  &\, (P_t^{D,n} |f|^p(x))^{\ff 1 p} \| (\si^n)^* (a^n)^{-1}\|_\infty C_{\bb,1,{t}}^{\ff 1 {2}}
\bigg(\E\sup_{s\in [0,t]} |\nn_v X_s^{n,x}|^{\ff{2p}{p-2}}\bigg)^{\ff{p-2}{2p}} \end{split}\end{equation}
for any $x\in D,\,v\in\R^d,\,t\in(0,T],\,f\in\B_b(\bar D)$ and $n\geq1$,  where  $\bb_s$ is an adapted real process with $\bb_0=1$ and $\bb_s=0$ for $s\ge t\land\tau_n$ will be constructed such that \eqref{BB2} holds for $\vv=1$.

(b) Let $f_i^n:= \int_{\R\times\R^d} f_i (x+y)g_n(s,y)\,\d s,\ 1\le i\le k_2.$ Then
$$\| \nn a^n\|\le \sum_{i=1}^{k_2} f_i^n, $$
and
\begin{equation}\label{capr}
\begin{aligned}
&\lim_{n\to\8}\left\{\|b^{1,n}-b^{(1)}\|_{\tt L^{p_1}_{q_1}(T)}+
\sum_{i=1}^{k_2}\|f^n_i-f_i\|_{\tt L^{p_i}_{q_i}(T)} \right\}=0,\\
&\lim_{n\to\8}\sup_{t\in[0,T]}\sup_{x\in\R^d}\|\si^n_t(x)-\si_t(x)\|_{HS}=0,\\
&\sup_{n\geq1}\left\{\|a^n\|_\infty+\|(a^n)^{-1}\|_\infty\right\}\leq \|a\|_\infty+\|a^{-1}\|_\infty,\\
&\sup_{n\geq1}\left\{\|b^{1,n}\|_{\tt L^{p_1}_{q_1} (T) }+ \sum_{i=1}^{k_2} \|f_i^n\|_{\tt L^{p_i'}_{q_i'}  (T)}\right\}\leq \|b^{(1)}\|_{\tt L^{p_1}_{q _1} (T)}+\sum_{i=1}^{k_2} \|f_i\|_{\tt L_{q_i'}^{p_i'} (T)}.
\end{aligned}
\end{equation}
By the proof of \cite[Theorem 3.9]{XZ}, we have
\begin{align}\label{xnx}
\lim_{n\to\8}\E\left[\sup_{t\in[0,T]}|X_t^{n,x}-X_t^x|\right]=0,\ \  x\in \R^d.
\end{align}
In addition, as shown in the proof of	\cite[Theorem 2.1(1)]{W23b}, for any $k\in [1,\infty)$, one   finds a constant $c(k)>0$ such that
$$\sup_{n\geq1}\E\left[\sup_{t\in[0,T]}|\nn_v X^{n,x}_t|^k\right] \le |v|^k c(k),\ \ \  x, v\in\R^d.$$
This together with \eqref{GRDn}   implies that for  for any  $p>2$, there exists  a constant $c_1>0$ uniformly in $n$ such that 
\beq\label{GRDn'}     \left|\nn  P_t^{D,n}f(x)\right|\le  c_1 (P_t^{D,n} |f|^p(x))^{\ff 1 p} \bigg(\E \int_0^t|\bb_s'|^2\d s \bigg)^{\ff 1 {2}} \end{equation}
holds for any   $t\in(0,T],\,x\in D,\,f\in\B_b(\bar D)$ and $ n\geq1.$ 

(c) To derive explicit gradient estimate from \eqref{GRDn'}, we construct the process $\bb$ by following the line of \cite{AW}.
Since $\pp D$ is uniformly of class $C^{2+\aa}$, we find a constant $r_1\in (0,r_0)$ such that $\rr_\pp \in C_b^2(D_{r_1})$, {see \cite[Proposition B.3]{FMP}}.
Let $h\in C_b^\infty([0,\infty))$ such that
\beq\label{GB} 0\le h\le 1,\ \ h'\ge 0,\ \ h(r)= r\ \text{for}\ r\le \ff{r_1}2,\ \ h(r)=1\ \text{for}\ r\ge r_1.\end{equation}
Then $g:=h\circ\rr_\pp\in C_b^2(\bar D)$ and there exists a constant $C_0>1$ such that
\beq\label{GC}C_0^{-1}  \le \ff{g}{\rr_\pp\land 1}(x)\le C_0,\ \ \ x\in \bar D.\end{equation}
By {\bf (C)}, $|b|+\|a\|+\|\nn^2\rr_\pp\|$ is bounded in $D_{r_0}$. Then there exists large enough $n_0\in \mathbb N$ such that
the family $\{|b^n|+\|a^n\|\}_{n\ge n_0}$ is uniformly  bounded on $[0,T]\times D_{r_1}$. So, there exists a constant $C_1>0$ such that
\beq\label{GC2} L_t^n g(x):= \ff 1 2{\rm tr}\{a_t^n(x)\nn^2g(x)\}+ b_t^n (x)\cdot\nn g(x)\le C_1,\ \ (t,x)\in [0,T]\times {D_{r_1}},\ n\ge n_0.\end{equation}
Now,   for any fixed $x\in D$, simply denote $X_t^{n,x}=X_t^n, \tau_n(x)=\tau_n$ and define
\beg{align*}
&S_n(t)=\int_{0}^{t}g^{-2}(X^n_{s\wedge\tau})\d s,\
&  s_n(t)=\inf\{s\in[0,T]: S_n(s)\geq t\},\ \ t\in[0,T]\ ,n\geq1.
\end{align*}
By $0\le g\le 1,$ \eqref{GC} and \eqref{GC2}, the same argument in \cite[Section 4]{AW} implies
that   $S_n(t)$ is increasing with $S_n(t)=\infty$ for $t>\tau_n$,  $s_n(t)\le t\land\tau_n,$
$S_n(s_n(t))=t$ for $t\leq S_n(\tau_n)$,    $s_n(S_n(t))=t$ for $t\leq \tau_n$.
Hence, for any   $t\in(0,T]$, 
$$
\beta_s:=1- \ff 1 t\int_0^{s\land s_n(t)}  g^{-2}(X_r^n) \d r,\ \  s\in [0,t]$$ satisfies
$\bb_0=1$,    and when $s\ge \tau_n \ (\ge s_n(t))$,
\begin{align*}
\beta_s=1-\ff 1 t\int_0^{s_n(t)}g^{-2}(X_r^n)\,\d r={1-}\ff 1 t S_n(s_n(t))={0}.
\end{align*}
Moreover,  because of    \eqref{GB}, $g:=h\circ\rr_\pp$ and \eqref{GC2},  by repeating  the calculations on page 118 in \cite{AW} for $(L_t^n, g)$ in pace of $(L,f)$ therein, we find a constant $c_3>0$ such that
$$\E\left[\int_0^{t\land \tau_n}|\beta_s'|^2\,\d s\right]=  \ff 1 {t^2}\E\left[\int_0^{s_n(t)}g^{-4}(X_s^n)\,\d s\right] \leq \ff{c_3}{t g^{2}(x)},\ \ \ n\ge n_0.$$
Combining this with \eqref{GRDn'} and \eqref{GC},     for any $p>2$ we find   a constant $c>0$ such that
\beq\label{GRDn''}     \left|\nn  P_t^{D,n}f(x)\right|\le  \ff c  {\ss t (\rr_\pp(x)\land 1)}(P_t^{D,n} |f|^p(x))^{\ff 1 p},\ \ n\ge n_0,\, f\in\B_b(\bar D).   \end{equation}

(d) Let $P_t^{D,n}f|_{D^c}=0$  for $t>0$, so that $ \nn P_t^{D,n}f|_{D^c}=0$. Then   \eqref{GRDn''} is equivalent to
\beg{align*}&\big|P_t^{D,n}f(x)-P_t^{D,n}f(y)\big|\\
&\le \ff{c|x-y| }  {\ss t }\int_0^1 \ff 1 { \rr_\pp(sx+(1-s)y)\land 1}\big(P_t^{D,n} |f|^p(sx+(1-s)y)\big)^{\ff 1 p}\d s
\end{align*}
for any  $n\geq n_0 ,\,t\in(0,T],\,x,y\in \bar D$ and $f\in\B_b(\bar D)$.  {Combining with the flow property  \eqref{FL}}, to deduce the desired gradient estimate on $P_t^D f$, it  remains to show
\beq\label{LMT} P_t^Df=\lim_{n\to\infty} P_t^{D,n} f,\ \ \ t\in (0,T],\ f\in C_b(\bar D).\end{equation}
For any $t\in(0,T],\,x\in D$   and $n\geq1$,
\begin{equation*}
\begin{aligned}
&\P\big(t<\tau,t\geq \tau_n \big)=\P\left( \inf_{s\in[0,t]}{\rho_\pp(X_s)}>0,\inf_{s\in[0,t]}{\rho_\pp(X_s^n)}=0\right) \\
&\leq \P\left(\sup_{s\in[0,t]}|X_s^{n}-X_s|\geq\vv\right)+\P\left(0<\inf_{s\in[0,t]}{\rho_\pp(X_s)}<\vv\right).
\end{aligned}
\end{equation*}
By first  letting $n\uparrow\8$ and then $\vv\downarrow0$, it follows from \eqref{xnx}  that \begin{align*}
\lim_{n\uparrow\8}	\P\big(t<\tau,t\geq \tau_n\big)=0,\ \ t\in(0,T],\,x\in D.
\end{align*}
Similarly,  by $\P(\tau =t)=0$ and \eqref{xnx},  {letting  $$\rr(x,D):=\inf_{y\in D}|x-y|,\ \ x\in \bar D,$$ }we have
\begin{equation*}
\begin{aligned}
&\P\big(t\geq\tau,t< \tau_n\big)=\P\big(t>\tau,t< \tau_n\big)\\
&=\P\left( \sup_{s\in[0,t]}\rr(X_s,D)>0,\sup_{s\in[0,t]}\rr(X_s^{n},D)=0\right) \\
&\leq \P\left(\sup_{s\in[0,t]}|X_s^{n}-X_s|\geq\vv\right)+\P\left(0<\sup_{s\in[0,t]}\rr(X_s,D)<\vv\right)
\end{aligned}
\end{equation*}
for any $t\in(0,T],\,x\in D$ and $n\geq1,$ which  tends to $0$ by letting first $n\uparrow\8$ then $\vv\downarrow0$. Hence,
\begin{align*}
\1_{\{t<\tau_n\}}\xrightarrow{\text{a.s.}}\1_{\{t<\tau\}} \text{ as }n\to\8,\ \ t\in(0,T],\,x\in D. 
\end{align*}
By combining this with   \eqref{xnx}  and the dominated convergence theorem, we derive     \eqref{LMT}. Therefore, the proof is finished.
\end{proof}

\section{Gradient estimate for  killed DDSDEs}

In this part, we establish the   estimates   \eqref{G'} and \eqref{R} for the killed DDSDE \eqref{E1} in the regular and singular cases respectively.

\subsection{Well-posedness of \eqref{E1}}

The well-posedness of \eqref{E1} has been derived in \cite{W23a}, where the noise may also depend on distribution, and coefficients are Lispchitz continuous in the distribution variable with respect to
the (truncated) $1$-Wasserstein  distance. Below we present a result ensuring the well-posedness for the case that the noise is distribution-free and the drift is Lipschitz continuous in the distribution variable
with respect to the total variation  distance. This situation is not included in \cite{W23a}, although the   proof is more or less standard by using the fixed point argument in distribution parameter.

For any $\mu\in C^w([0,T];\scr P_D)$, the set of all weakly continuous maps from $[0,T]$ to $\scr P_D$, let  $a:=\si\si^*$ and 
$$b^\mu_t(x):= b_t(x,\mu_t),\ \ \ (t,x)\in [0,T]\times \R^d.$$

\emph{\beg{enumerate} \item[{\bf (D)}]      
There exists   $0\le \psi \in L^2([0,T])$ such that
\beq\label{TV} |b_t(x,\mu)-b_t(x,\nu)|\le \psi_t \|\mu-\nu\|_{var},\ \ \ t\in [0,T],\   x\in \R^d, \ \mu,\nu\in\scr P_D.\end{equation} Moreover, one of the following conditions holds for any $\mu\in  C^w([0,T];\scr P_D):$  
\item[$\bullet$]  {\bf (A)} holds for $b^\mu$ in place of $b$ with constant $K$ possibly depend on $\mu$, or 
\item[$\bullet$]  $\|a\|_\8+\|a^{-1}\|_\8<\8$, and $\si_t$ and $b_t^\mu$  are continuous on $\R^d$ such that
for some $0\le h\in L^1([0,T])$ possibly depend on $\mu$, 
$$\<b_t^\mu(x)-b_t^\mu(y), x-y\>+\ff 1 2 \|\si_t(x)-\si_t(y)\|_{HS}^2\le h_t|x-y|^2,\ \ t\in [0,T],\ x,y\in\R^d.$$ 
\end{enumerate} }

\beg{thm}\label{TA}  Assume {\bf (D)}. Then the SDE $\eqref{E1}$ is well-posed, and there exists a constant $c>0$ such that the associated map $P_t^{D*}: \scr P_D\to \scr P_D$ satisfies
\beq\label{RG} \|P_t^{D*}\gg-P_t^{D*}\tt\gg \|_{var}\le c\|\gg-\tt\gg\|_{var},\ \ \ t\in [0,T],\ \gg,\tt\gg\in \scr P_D.\end{equation} \end{thm}

\beg{proof}
(a) Let $X_0$ be a fixed $\F_0$-measurable random variable on $\bar D$ with $\gg:=\L_{X_0}^D\in \scr P_D.$ Under {\bf (D)}, for any
$$\mu\in \scr C_T^\gg:=\big\{\mu\in C^w([0,T]; \scr P_D):\ \mu_0=\gg\big\},$$ the SDE
\beq\label{DCS} \d X_t^\mu= b_t^\mu(X_t^\mu)\d t+  \si_t  (X_t^\mu)\d W_t,\ \ X_0^\mu=X_0,\ t\in [0,T]\end{equation}  is well-posed,
which follows from \cite[Proposition 5.2]{HRW} when {\bf (A)} holds for $b^\mu$ in place of $b$, and is well-known in the other monotone case.
Let
$$\Phi_t\mu:=\L^D_{X_t^\mu},\ \ t\in [0,T].$$
Then $\Phi: \C_T^\gg \to \C_T^\gg,$
and   for the well-posedness it suffices to show that  $\Phi$ has a unique fixed point $\bar\mu\in \C_T^\gg$, and in this case $X_t:= X_{t\land \tau(X^{\bar\mu})}^{\bar \mu}$ becomes the unique solution of \eqref{E1},
where
$$ \tau(X^{\bar\mu}):=\inf\{t\ge 0:\ X_t^{\bar\mu}\in\pp D\}. $$
To this end, we will find $\ll>0$ such that $\Phi$ is contractive under the following complete metric on $\C_T^\gg$:
$$\| \mu- \nu\|_{\ll}:= \sup_{t\in [0,T]} \e^{-\ll t} \| \mu_t-\nu_t\|_{var},\ \ \  \mu,\nu\in \C_T^\gg.$$

For a different element $\nu\in \C_T^\gg,$ we reformulate \eqref{DCS} as
$$ \d X_t^\mu= b_t^\nu(X_t^\mu)\d t+  \si_t (X_t^\mu)\d \tt W_t,\ \ X_0^\mu=X_0,\ t\in [0,T],$$
where, by \eqref{TV} and Girsanov's theorem,
$$\tt W_t:= W_t-\int_0^t (\si_s^*a_s^{-1})(X_s^\mu) \big\{b_s^\nu(X_s^\mu)- b_s^\mu(X_s^\mu)\big\}\d s,\ \ t\in [0,T]$$
is an $m$-dimensional Brownian motion under the weighted probability
$  \d \Q:= R_T\d\P,$ for the exponential martingale
$$R_t:=\e^{\int_0^t\<\xi_s,\d W_s\>- \ff 12\int_0^t |\xi_s|^2\d s},\ \ \ \xi_s:= (\si_s^*a_s^{-1})(X_s^\mu) \big\{b^\nu_s (X_s^\mu)- b_s^\mu(X_s^\mu)\big\},\ \ \ 0\le s\le t\le T.$$
So, by the weak uniqueness of \eqref{DCS},
$$\|\Phi_t\mu-\Phi_t\nu\|_{var}=  \sup_{\|f\|_\infty\le 1} \big|\E\big[1_{\{t<\tau(X^\mu)\}} (R_t-1) f(X_t^\mu)\big]\big|\le \E[|R_t-1|],\ \  t\in[0,T]. $$
This together with Pinsker's inequality and \eqref{TV} implies that
\beq\label{PV1} \|\Phi_t\mu-\Phi_t\nu\|_{var} \le \ss{2\E[R_t\log R_t]}\le \bigg( c_0 \int_0^t \psi_s^2 \|\mu_s-\nu_s\|_{var}^2\d s\bigg)^{\ff 1 2},\ \ t\in [0,T]\end{equation} holds for some constant $c_0>0$. 
Thus, when $\ll>0$ is large enough,
\beg{align*}&\|\Phi\mu-\Phi\nu\|_{\ll}:= \sup_{t\in [0,T]} \e^{-\ll t} \|\Phi_t\mu-\Phi_t\nu\|_{var} \\
&\le  \|\mu-\nu\|_\ll \sup_{t\in [0,T] } \bigg( c_0 {\int_0^t} \psi_s^2\e^{-2(t-s)}\d s\bigg)^{\ff 1 2}\le \ff 1 2 \|\mu-\nu\|_\ll,\ \ \mu,\nu\in \C_T^\gg.\end{align*}
Therefore,  $\Phi$ has a unique fixed point in $\C_T^\gg.$

(b) Let $\tt\gg\in \scr P_D$ be different from $\gg$, and choose
$$\mu_t=P_t^{ D*}\gg,\ \ \ \nu_t=P_t^{ D*}\tt\gg,\ \ \ t\in [0,T].$$
Let  $X_0$ and $Y_0$ be $\F_0$-measurable random variables   on $\bar D$ such that
$$\L_{X_0}^D=\gg,\ \ \L_{Y_0}^D = \tt \gg,\ \ \ \P\big(\1_{\{X_0\in D\}}X_0\ne \1_{\{Y_0\in D\}}Y_0\big) \le 2\|\gg-\tt\gg\|_{var}.$$
Let $\tt X_t^\mu$ and  $\tt X_t^\nu$, with $\tt X_0^\mu=X_0$ and  $\tt X_0^\nu=Y_0$,  solve the following SDE
\beq\label{DCSS} \d \tt X_t = b_t^\nu(\tt X_t)\d t+\si_t(\tt X_t )\d W_t,\ \ \ t\in [0,T].\end{equation}
Then
$$P_t^{D*}\gg= \Phi_t\mu:=\L^D_{X_t^\mu},\ \ \ \Phi_t\nu= \L_{\tt X_t^\mu}^D,\ \ \  P_t^{D*}\tt\gg= \L^D_{\tt X_t^\nu},\ \ \ t\in [0,T].$$
Moreover,    the pathwise uniqueness of \eqref{DCSS}  implies
$$\big\{1_{\{  X_0\in D\}}X_0 = \1_{\{Y_0\in D\}}Y_0\big\}\subset \big\{1_{\{t<\tau(\tt X^\mu)\}}\tt X_t^\mu= 1_{\{t<\tau(\tt X^\nu)\}}\tt X_t^\nu\big\},$$
so that
$$\|\Phi_t\nu- P_t^{D*}\tt\gg\|_{var}= \|\L_{\tt X_t^\mu}^D-\L_{\tt X_t^\nu}^D\|_{var}\le \|\gg-\tt\gg\|_{var}.$$
Combining this with \eqref{PV1}, we obtain
\beg{align*}& \|P_t^{D*}\gg-P_t^{D*}\tt\gg\|_{var} \le  \|\Phi_t\mu-\Phi_t\nu\|_{var} +  \|\Phi_t\nu- P_t^{D*}\tt\gg \|_{var} \\
&\le   \bigg( c_0 \int_0^t \psi_s^2 \|P_s^{D*}\gg -P_s^{D*}\tt\gg\|_{var}^2\d s\bigg)^{\ff 1 2} + \|\gg-\tt\gg\|_{var},\ \ t\in [0,T].\end{align*}
Squaring both sides and applying Gronwall's inequality, we derive  \eqref{RG} for some constant $c>0.$

\end{proof}

\subsection{Gradient estimate  in the regular case}

\beg{thm} \label{TA2}
If {\bf (B)} holds for $b^\mu$ in place of $b$ uniformly in $\mu\in  C^w([0,T];\scr P_D),$ and there exists a constant $K>0$ such that
\beq\label{TV'} |b_t(x,\mu)-b_t(x,\nu)|\le K \|\mu-\nu\|_{var},\ \ \ t\in [0,T],\   x\in\R^d, \ \mu,\nu\in\scr P_D.\end{equation}  Then  there exists a constant $c>0$ such that  the   map $P_t^{D*}: \scr P_D\to \scr P_D$
associated with $\eqref{E1}$ satisfies
\beq\label{RR} \|P_t^{D*}\mu_0-P_t^{D*}\nu_0\|_{var}\le \ff{c}{ \ss t} \hat\W_1(\mu_0,\nu_0),\ \ \ t\in (0,T],\ \mu_0,\nu_0\in \scr P_D.\end{equation}  \end{thm}

\beg{proof} For $\mu_0, \nu_0 \in \scr P_D$, define
$$\mu_t:=P_t^{D*}\mu_0,\ \ \ \nu_t:= P_t^{D*}\nu_0,\ \ \ t\in [0,T],$$
and let
$P_{s,t}^{D,\mu}$  (respectively  $P_{s,t}^{D,\nu}$) be the Dirichlet semigroup $P_{s,t}^D$ defined in \eqref{TS} for $b^\mu$ (respectively $b^\nu$) in place of $b$.
By Lemma \ref{LM}, there exists a constant $c_1>0$ such that
\beq\label{KB1} \|\nn P_{s,t}^{D,\mu}  f\|_\infty\le \ff {c_1\|f\|_\infty}{\ss {t-s}},\ \ \ 0\le s<t\le T,\ f\in \B_b(\bar D),\end{equation}
and for any $t\in(0,T]$, $f\in\B_b(\bar D)$, 
\beq\label{BK2} \pp_s  P_{s,t}^{D,\mu}  f= - \Big\{\ff 1 2 {\rm tr}(a_s\nn^2) + b_s^\mu \cdot\nn\Big\} P_{s,t}^{D,\mu}  f,\ \ \ s\in [0, t).\end{equation}
For any $\gg\in \scr P_D$ and $t\in [0,T],$ let $P_t^{D,\mu*}\gg\in \scr P_D$ be defined by
$$(P_t^{D,\mu*}\gg)(A):= \gamma\big( P_{t}^{D,\mu} 1_A\big),\ \ \ A\in\B_D,  $$ where $ P_{t}^{D,\mu} := P_{0,t}^{D,\mu}$.
By the same reason deducing \eqref{R} from \eqref{*G},     \eqref{KB1} implies
\beq\label{KB3} \|P_t^{D,\mu*}\mu_0- P_t^{D,\mu*}\nu_0\|_{var}\le \ff {c_2}{\ss {t}}\hat\W_1(\mu_0,\nu_0), \ \ t\in (0,T]\end{equation}
for some constant $c_2>0$ uniformly in $\mu_0,\nu_0.$
On the other hand, let $X_s^\nu$ solve the SDE
$$\d X_s^\nu= b_s^\nu(X_s^\nu)\d s+\si_s(X_s^\nu)\d W_s,\ \  \L_{X_0^\nu}^D=\nu_0,\  s \in [0,T].$$
For any  $f\in \B_b(\bar D)$, by \eqref{BK2} and It\^o's formula, we obtain
$$\d  P_{s,t}^{D,\mu}f(X_s^{\nu})= \big\<b_s^\nu(X_s^{\nu})- b_s^\mu(X_s^{\nu}), \nn  P_{s,t}^{D,\mu} f(X_s^{\nu})\big\>\d s+\d M_s,\ \  0\leq s<t\leq T,$$
where by \eqref{KB1} and $\|\si\|_\infty<\infty$,
$$\d M_s:= \big\<\nn P_{s,t}^{D,\mu}  f(X_s^{\nu}), \si_s(X_s^{\nu})\d W_s\big\>,\ \
0\leq s<t\leq T $$ is a martingale.
Combining this with
$$P_{s\land \tau(X^{\nu}),t}^{D,\mu }f(X_{s\land \tau(X^\nu)}^{\nu}) =0\ \text{if}\ t>s\ge \tau(X^\nu),$$
we derive the Duhamel formula
\beg{align*}& \nu_0(P_t^{D,\nu}f) =\lim_{s\uparrow t} \E\big[\1_{\{s<\tau(X^\nu)\}}P_{s,t}^{D,\mu } f(X_s^{\nu})\big] = \lim_{s\uparrow t} \E\big[P_{s\land \tau(X^\nu),t}^{D,\mu} f(X_{s\land \tau(X^\nu)}^{\nu})\big] \\
&= \E[P_{0,t}^{D,\mu } f(X_0^{\nu})]+\E \int_0^{t\land\tau(X^\nu)}  \big\<b_s^\nu- b_s^\mu, \nn P_{s,t}^{D,\mu } f \big\>(X_s^\nu)\d s  \\
& = \nu_0 ( P_t^{D,\mu } f)+\int_0^t \nu_0 \left(  P_s^{D,\nu}\big\<b_s^\nu- b_s^\mu, \nn P_{s,t}^{D,\mu } f  \big\>\right)\d s,\ \ \ t\in [0,T],\, f\in\B_b(\bar D). \end{align*}
Combining this with \eqref{TV'},   \eqref{KB1} 
and noting that
$$\mu_t:=P_t^{D*}\mu_0= P_t^{D,\mu*}\mu_0,\ \ \ \nu_t:=P_t^{D*}\nu_0= P_t^{D,\nu*}\nu_0,\ \  t\in[0,T], $$ we obtain
$$ \|P_t^{D,\mu*}\nu_0- \nu_t\|_{var}=\sup_{\|f\|_\infty\le 1} \big|\nu_0(P_t^{D,\mu } f-P_t^{D,\nu }f)\big| \le c_2K \int_0^t \ff{\|\mu_s-\nu_s\|_{var} }{\ss{t-s}}\d s,\ \ t\in [0,T].$$
Combining this with \eqref{KB3} we arrive at
\beq\label{KB4}\beg{split} \|\mu_t- \nu_t\|_{var}& \le \| \nu_t- P_t^{D,\mu*}\nu_0\|_{var}+ \|P_t^{D,\mu*}\nu_0- P_t^{D,\mu*}\mu_0\|_{var}\\
&\le c_2K \int_0^t \ff{\|\mu_s-\nu_s\|_{var} }{\ss{t-s}}\d s+ \ff {c_2}{\ss {t}}\hat\W_1(\mu_0,\nu_0),\ \ \ t\in (0,T].\end{split} \end{equation}
It is trivial that
$$\phi_\ll:=\sup_{t\in [0,T]} \e^{-\ll t} \ss t \|\mu_t- \nu_t\|_{var} <\infty,\ \ \ \ll\in (0,\infty).$$
By \eqref{KB4}, when $\ll>0$ is large enough,
$$\phi_\ll \le c_2 \hat\W_1(\mu_0,\nu_0)+  c_2K \phi_\ll \sup_{t\in (0,T]}  \int_0^t \ff{\e^{-\ll(t-s)} \ss t}{\ss{s(t-s)}}\d s\le c_2 \hat\W_1(\mu_0,\nu_0)+  \ff 1 2 \phi_\ll,$$
so that
$\phi_\ll\le 2  c_2 \hat\W_1(\mu_0,\nu_0).$ Consequently,
$$\|P_t^{D*}\mu_0- P_t^{D*}\nu_0\|_{var}=\|\mu_t- \nu_t\|_{var}\le \ff{2  c_2\e^{\ll T}} {\ss t}\hat\W_1(\mu_0,\nu_0),\ \ t\in (0,T],$$
so that \eqref{RR} holds for $c:= 2  c_2\e^{\ll T}.$
\end{proof}

\subsection{Gradient estimate   in the singular case}

We make the following assumption.
\emph{\beg{enumerate} \item[{\bf (E)}]    Assume that {\bf (C)} holds uniformly in $\mu\in \C^w([0,T];\scr P_D)$, and there exist $k_3 \in \mathbb N$ and functions $\{0\le \tt f_i\in \tt L_{\tt q_i}^{\tt p_i}\big\}$ for some
$\{(\tt p_i,\tt q_i):\ 1\le i\le k_3\big\}\subset \scr K$, such that
\beq\label{Q*}  |b_t(x,\mu)- b_t(x,\nu)|\le \hat\W_1(\mu,\nu) \sum_{i=1}^{k_3} \tt f_i(t,x),\ \ (t,x,\mu, \nu)\in [0,T]\times \bar D\times \scr P_D\times \scr P_D.\end{equation}
\end{enumerate}}
\beg{thm}  Assume {\bf (E)}.
Then  there exists a constant $c>0$ such that  the   map $P_t^{D*}: \scr P_D\to \scr P_D$
associated with $\eqref{E1}$ satisfies
\beq\label{GG} \|P_t^{D*}\mu_0-P_t^{D*}\nu_0\|_{var}\le \ff{c}{ \ss t} \hat\W_1^\rr(\mu_0,\nu_0),\ \ \ t\in (0,T],\ \mu_0,\nu_0\in \scr P_D.\end{equation}
\end{thm}

\beg{proof}  For any $\mu_0,\nu_0\in \scr P_D$,  let $\mu_t,\nu_t, P_t^{\mu*} $ and $P_t^{\nu*}$  be in the proof of Theorem \ref{TA2}.

By {\cite[Theorem 3.1]{W23a}}, which is proved for $k_1=1$, but the proof also works for $k_1\in \mathbb N$  due to \cite[Proposition 5.2]{HRW}, {the SDE \eqref{E1} is well-posed, and }
we can find a constant $c_0>0$ such that
\beq\label{X0} \hat\W_1(\mu_t,\nu_t) \le c_0\hat\W_1(\mu_0,\nu_0)  \le c_0\hat\W_1^\rr(\mu_0,\nu_0),\ \ \ t\in [0,T].\end{equation}
By Theorem \ref{T2} for $b^\mu$ in place of $b$, so that \eqref{GF} holds for $P_t^{D,\mu }$ in place of $P_t^D$, hence  there exists a constant $c_1>0$ uniformly in $x$ 
such that
$$|\nn P_t^{D,\mu } f(x)|  \le \ff{c_1\|f\|_\infty}{\ss t (1\land\rr_\pp(x))},\ \ t\in (0,T],\ x\in D,\ f\in \B_b(\bar D).$$
By the same reason deducing  \eqref{G'} from \eqref{G}, this implies
\beq\label{X2}  \|P_t^{D,\mu*} \mu_0- P_t^{D,\mu*}\nu_0\|_{var} \le \ff{c_1}{\ss t} \hat \W_1^\rr(\mu_0,\nu_0),\ \ \ t\in (0,T].\end{equation}
Moreover,  by  \eqref{Q*},  \eqref{X0}  and the Girsanov theorem  as in the proof of  \eqref{PV1}, we find constants $c_2,c_3>0$ such that
\beq\label{XX}
\begin{aligned}
\|P_t^{D,\mu*}\nu_0-P_t^{D,\nu*}\nu_0\|_{var}\le& c_2 \bigg(\sum_{i=1}^{ k_3} \E_\Q \int_0^t  \tt f_i(s, X_s^\mu)^2  \hat \W_1(\mu_s,\nu_s)^2\d s\bigg)^{\ff 1 2}\\
\le& c_3\hat\W^1(\mu_0,\nu_0)\bigg(\sum_{i=1}^{k_3} \E_\Q \int_0^t  \tt f_i(s, X_s^\mu)^2  \d s\bigg)^{\ff 1 2},\ \ t\in [0,T],
\end{aligned} \end{equation}
where $X_s^\mu$ solves \eqref{DCS} for $\L_{X_0^\mu}^D=\nu_0$, and $\d\Q=R_T\d\P$ as defined before \eqref{PV1}. Noting that
the  law of $X^\mu$ under $\Q$ coincides with that of $X^\nu$ under $\P$, where $X_s^\nu$ solves \eqref{DCS} with $\L_{X_0^\nu}^D=\nu_0$, we have
$$\E_\Q \int_0^t  \tt f_i(s, X_s^\mu)^2\d s= \E  \int_0^t  \tt f_i(s, X_s^\nu)^2\d s,$$
so that by Krylov's estimate,   see \cite[Lemma 1.2.3]{RW24}, we find a constant $c_4>0$ such that
$$\sum_{i=1}^k  \E_\Q \int_0^t  \tt f_i(s, X_s^\mu)^2\d s\le c_4.$$
Combining this with  \eqref{XX}, we find a constant $c_5>0$ such that
$$\|P_t^{D,\mu*}\nu_0-P_t^{D,\nu*}\nu_0\|_{var}\le c_5 \hat \W_1(\mu_0,\nu_0),\ \ t\in [0,T].$$
This together with \eqref{X2}  and the fact that $\hat\W_1\leq\hat\W^\rr_1$  proves \eqref{GG} for some constant $c>0$,  since
\beg{align*}&\|P_t^{D*}\mu_0-P_t^{D*}\nu_0\|_{var}= \|\mu_t-\nu_t\|_{var} =\|P_t^{D,\mu*} \mu_0- P_t^{D,\nu*}\nu_0\|_{var}\\
&\le \|P_t^{D,\mu*} \mu_0- P_t^{D,\mu*}\nu_0\|_{var}+  \|P_t^{D,\mu*}\nu_0-P_t^{D,\nu*}\nu_0\|_{var}, \ \  t\in[0,T].\end{align*}

\end{proof}

\paragraph{Acknowledgement.} The authors would like to thank the referees for corrections and valuable comments.  
\section*{Declarations}
\noindent\textbf{Funding:} This work is supported in part by the National Key R\&D Program of China (No. 2022YFA1006000), and State Key Laboratory of Synthetic Biology, Tianjin University.

\

\noindent\textbf{Conflict of interest:} The authors declare that they have no known competing interests.

\end{document}